\documentclass{article}

\usepackage[dvipsnames]{xcolor}
  \definecolor{myurlcolor}{rgb}{0.6,0,0}
  \definecolor{mycitecolor}{rgb}{0,0,0.8}
  \definecolor{myrefcolor}{rgb}{0,0,0.8}

\usepackage{amsfonts}
\usepackage{amssymb}  
\usepackage{amsthm} 
\usepackage{amsmath}
\usepackage{authblk}
\usepackage[font=footnotesize]{caption}
\usepackage{circuitikz}
\usepackage[inline]{enumitem}
  \setlist{itemsep=0em, topsep=0em, parsep=0em}
\usepackage{etoolbox}
\usepackage{hyperref}
\hypersetup{colorlinks,
	linkcolor=myrefcolor,
	citecolor=mycitecolor,
	urlcolor=myurlcolor}
\usepackage{graphicx}
  \graphicspath{ {assets/} }
  \usepackage{mathtools}
\usepackage{natbib}
\usepackage{stmaryrd} 
\usepackage{tikz}
  \usetikzlibrary{matrix,
    arrows,
    shapes,
    decorations.markings,
    decorations.pathreplacing}
\usepackage{todonotes}

\renewcommand{\epsilon}{\varepsilon}
\newcommand{\op}{^{\scriptsize{ \textrm{op} } }}

\newcommand{\bydef}{\coloneqq}

\newcommand{\NN}{\mathbb{N}}

\newcommand{\A}{\cat{A}}

\newcommand{\C}{\cat{C}}

\newcommand{\X}{\cat{X}}

\newcommand{\Set}{\cat{Set}}

\newcommand{\Graph}{\cat{Graph}}

\newcommand{\HGraph}{\cat{HGraph}}
\newcommand{\DHGraph}{\cat{DHGraph}}

\newcommand{\Topos}{\cat{Topos}}

\newcommand{\Csp}{\cat{Csp}}

\newcommand{\StrCsp}{\cat{StrCsp}}

\newcommand{\Lang}{\mathrm{Lang}}
\newcommand{\LLang}{\bicat{Lang}}
\newcommand{\Comp}{\mathrm{Comp}}

\newcommand{\df}[1]{\textbf{#1}}
\newcommand{\cat}[1]{\mathsf{#1}}
\newcommand{\bicat}[1]{\mathbf{#1}}

\newcommand{\type}[1]{\mathtt{#1}}

\newcommand{\from}{\colon}

\newcommand{\xto}[1]{\xrightarrow{#1}}

\newcommand{\dderiv}[2]{#1 \rightsquigarrow #2}
\newcommand{\deriv}[2]{#1 \rightsquigarrow^\ast #2}
\renewcommand{\gets}{\leftarrow}

\newcommand{\xgets}[1]{\xleftarrow{#1}}
\newcommand{\spn}[3]{#1 \gets #2 \to #3}

\newcommand{\csp}[3]{#1 \to #2 \gets #3}

\DeclareMathOperator{\id}{id}

\newtheorem{therm}{Theorem}[section]
\newtheorem{lemma}[therm]{Lemma}
\newtheorem{proposition}[therm]{Proposition}
\newtheorem{corollary}[therm]{Corollary}

\theoremstyle{remark}

\theoremstyle{definition}
\newtheorem{definition}[therm]{Definition}

\newtheorem{example}[therm]{Example}

\tikzset{
  cd/.style={
    ->,
    scale=6,
    >=angle 90,
    font=\scriptsize}}

\tikzset{
  graph/.style={
    ->,
    -{Stealth[scale=0.75]},
    font=\scriptsize}}

\tikzset{->-/.style={decoration={%
      markings,
      mark=at position .5 with {\arrow{>}}},
      postaction={decorate}}}

\tikzset{->-pos/.style={decoration={%
      markings,
      mark=at position #1 with {\arrow{>}}},
      postaction={decorate}}}

\tikzset{-|->/.style={decoration={%
      markings,
      mark=at position .5 with {\arrow{|}},
      mark=at position 1 with {\arrow{>}}},
    postaction={decorate}}}

\newlength\mylen
\newcommand{\horarrow}{%
  \to\kern-0.55\mylen\vline height 1.2ex depth
  -0.4pt\kern0.55\mylen}

\newcommand{\adjunction}[4]{%
  \begin{tikzpicture}[baseline=-3pt]
    \node (1) at (0,0) {\( #1 \)};
    \node (2) at (2,0) {\( #4 \)};
    \draw [->]
    ([yshift= 4pt]2.west) to
    node [above] {\scriptsize{ $ #2 $ }}
    ([yshift= 4pt]1.east);
    \draw [->]
    ([yshift= -4pt]1.east) to
    node [below] {\scriptsize{ $ #3 $ }}
    node [above,yshift= -1.5pt] {\scriptsize{$ \perp $}}
    ([yshift= -4pt]2.west);
  \end{tikzpicture}
}

\title{Rewriting Structured Cospans}
\author{Daniel Cicala}
\date{%
  Department of Mathematics and Statistics\\%
  Southern Connecticut State University\\%
  cicalad1@southernct.edu\\[2ex]%
  \today
} 

\begin{document}

\maketitle

\begin{abstract}
We develop a theory of rewriting for structured cospans in order to extend compositional methods for modeling open networks. First, we introduce a category whose objects are structured cospans, and establish conditions under which it is adhesive or a topos. These results guarantee that double pushout rewriting can be applied in this setting. We then define structured cospan grammars and construct their associated languages via a 2-categorical framework, capturing both network composition and rewrite dynamics. As an application, we show that for graphs, hypergraphs, Petri nets, and their typed variants, any grammar induces the same language as its corresponding discrete grammar. This equivalence enables an inductive characterization of rewriting, thereby generalizing classical results from graph transformation to a broader class of categorical models.
\end{abstract}

\section{Introduction} \label{sec:Intro}

Structured cospans are a framework for reasoning about networks with inputs and outputs. Rewriting is a topic that covers methods for editing substructures of an object such as a string or a graph. In this paper, we introduce rewriting to structured cospans.

Baez and Courser introduced structured cospans as an abstract framework to model open networks \cite{StructuredCospans, Courser_Open}. A structured cospan is a diagram of the form
\begin{equation}\label{eq:Intro_StrCspForm}
  La \to x \gets Lb
\end{equation}
where $ L \from \A \to \X $ is a functor whose codomain has chosen pushouts. This functor is a convenient bookkeeping device allowing us to separate open networks into the interface types $\A$ and the network types $ \X $.  We interpret \eqref{eq:Intro_StrCspForm} as a network $ x $ with inputs $ La $ and outputs $ Lb $. This network can then \emph{connect} to a compatible network, that is one with inputs $Lb$, by pushout: \[(La \to x \gets Lb \to y \gets Lc) \mapsto (La \to x +_{Lb}y \gets Lc).\] We then get a category $_L\Csp$ whose objects are those of $ \A $ and whose arrows are structured cospans $ La \to x \gets Lb $ where composition models the connecting of networks.

Given the ability to model open networks using structured cospans, we would like tools to analyze these models. In this paper, we adapt for structured cospans the tool of double pushout (DPO) rewriting, an algorithmic technique for creating a new object from an old object according to a given set of rules.

\subsection{Double Pushout Rewriting}
\label{sec:intro-dpo-rewriting}

Rewriting offers a method of replacement, often for the purpose of generating or simplifying. There are numerous examples.  In the case of electrical circuits, rewriting provides a method to replace two resistors in series with a single resistor. In database schemas, rewriting offers a way to transform one schema into another \cite{rodgers-graph-rewriting-database}. Rewriting generates strings in formal languages \cite{baader-term-rewriting}. 

Double pushout (DPO) rewriting is an algebraic mechanism that can perform the rewriting. It works in any adhesive category \cite{LackSobo_Adhesive}. One starts with a \emph{grammar} $ ( \C,P ) $, that is an adhesive category $ \C $ and a set $ P \bydef \{ \ell_j \gets k_j \to r_j \} $ of spans in $ \C $ with monic arrows. Each span is called a \emph{rewrite rule}, or just \emph{rule} for short. We interpret a rule $ \ell \gets k \to r $ as stating that $ r $ replaces $ \ell $ in a manner that fixes $ k $. This rule can be applied to any object $ \ell' $ by realizing a \emph{double pushout diagram}
\begin{equation} \label{eq:Intro_DpoDiagram}
 \begin{tikzpicture}[
  baseline=(current bounding box.center)]
  scale=0.75]
  \node (1) at (0,1) {$ \ell $};
  \node (2) at (1,1) {$ k $};
  \node (3) at (2,1) {$ r $};
  \node (4) at (0,0) {$ \ell' $};
  \node (5) at (1,0) {$ k' $};
  \node (6) at (2,0) {$ r' $};
  \draw [cd]
  (2) edge (1)
  (2) edge (3)
  (5) edge (4)
  (5) edge (6)
  (1) edge node[left]{$ m $} (4)
  (2) edge (5)
  (3) edge (6); 
  \draw (0.3,0.4) -- (0.4,0.4) -- (0.4,0.3);
  \draw (1.7,0.4) -- (1.6,0.4) -- (1.6,0.3);
\end{tikzpicture}
\end{equation}
where $ m $ identifies an instance of $ \ell $ in $ \ell' $ and  replaces it with $ r $. The resulting object of this rewriting process is $ r' $. The squares being pushouts means that the \emph{gluing condition} is met, which ensures the process outputs a well-defined object. For instance, in the case of graphs, the gluing condition means that there no node is removed leaving an edge unanchored. By considering all possible double pushout diagrams \eqref{eq:Intro_DpoDiagram} in $\C$ where the top row is in $P$, we get a relation $ \dderiv{\ell'}{r'}$ on the objects of $\C$. The primary item to study is the \emph{rewrite relation} $\deriv{}{}$ which is the reflexive and transitive closure of $\dderiv{}{}$. One may alternatively study the \emph{language} of a grammar, which is the category whose objects are those of $\C$ and arrows are generated by the rewrite relation. Both the rewrite relation and the language encode all the possible ways to rewrite one object into another.

\subsection{Rewriting Structured Cospans}
\label{sec:intro-rewriting-str-csps}

Observe that rewriting applies to the \emph{objects} of an adhesive category. So instead of working with the category $_L\Csp$ where structured cospans are arrows, we need to construct a category where they are objects.  We denote this by $ _L\StrCsp $. To build it from scratch, start with an adjunction $L \dashv R \from X \to A $.  We may consider that $\X$ contains every network type, $\A$ every interface type, $L$ includes each interface type into $\X$ as a sort of trivial network type, and $R$ returns what may loosely be considered the largest possible interface of a network.  Now, define $ _L\StrCsp $ to have  structured cospans for objects and commuting diagrams
\[
  \begin{tikzpicture}
    \node (1) at (0,1) {\( La \)};
    \node (2) at (1,1) {\( x \)};
    \node (3) at (2,1) {\( Lb \)};
    \node (4) at (0,0) {\( La'\)};
    \node (5) at (1,0) {\( x' \)};
    \node (6) at (2,0) {\( Lb' \)};
    \draw [cd]
    (1) edge node[]{} (2)
    (3) edge node[]{} (2)
    (4) edge node[]{} (5)
    (6) edge node[]{} (5)
    (1) edge node[left]{$Lf$}  (4)
    (2) edge node[left]{$g$}   (5)
    (3) edge node[right]{$Lh$} (6);
  \end{tikzpicture}
\]
as arrows. We show that if $\X$ and $\A$ are adhesive and $L$ preserves pullbacks, then $_L \StrCsp$ is adhesive (Theorem \ref{thm:strcsp-is-adhes}). Also, if we strengthen this to make $ L \dashv R $ a geometric morphism between topoi, then $ _L\StrCsp $ is a topos (Theorem \ref{thm:strcsp-istopos}).  It follows that, under modest conditions, $_L\StrCsp$ supports rewriting.

Though we restrict $ L $ more than Baez and Courser, our definition still covers many important examples.  One such example is the discrete graph functor $ L \from \Set \to \Graph $. More examples come by using slice categories $ \Graph / g $ for some graph $ g $ chosen to endow nodes and arrows with types. This was done to model the ZX-calculus \cite{ZX} and can be done to model passive linear circuits \cite{PassiveNets} by choosing $g$ to have a single node and edge set of possible resistances $(0,\infty)$. Other examples include directed hypergraphs, undirected hypergraphs, their typed versions, Petri nets, and marked Petri nets. 

To rewrite structured cospans, we begin with a \emph{structured cospan grammar} $ (_L\StrCsp, P) $, which is different than a grammar as discussed above because we require $ P $ to contain spans in $ _L\StrCsp $ of the form
\[
  \begin{tikzpicture}[baseline=(current bounding box.center)]
    \node (x') at (1,2) {$ x' $};
    \node (La) at (0,1) {$ La $};
    \node (x) at (1,1) {$ x $};
    \node (Lb) at (2,1) {$ Lb $};
    \node (x'') at (1,0) {$ x'' $};
    \draw [>->] (x) to (x');
    \draw [>->] (x) to (x'');
    \draw[cd]
    (La) edge[] (x)
    (Lb) edge[] (x)
    (La) edge[] (x')
    (Lb) edge[] (x')
    (La) edge[] (x'')
    (Lb) edge[] (x'');
  \end{tikzpicture}
\]
which rewrites $\csp{La}{x'}{Lb}$ into $\csp{La}{x''}{Lb}$. This condition is stronger than simply requiring monic-legged spans because the interfaces are fixed. 

Earlier, we associated to a grammar a 1-category called the ``language''. Now, with a structured cospan grammar $(_{L}\StrCsp,P)$, we construct its language $\LLang (_{L}\StrCsp,P)$ as a free 2-category on a computad \cite{Street_Computad}.  The reason for a 2-category is that there are two compositional structures at play. There is the composing of structured cospans and of the rewrite rules. These will respectively form the 1 and 2-morphisms of $\LLang (_{L}\StrCsp,P)$. 

\subsection{An Application: An Inductive Viewpoint on Rewriting}
\label{sec:intro-application}

In the classical topics of rewriting---formal languages and term rewriting---there are two approaches to defining the rewrite relation for a grammar. The first is an operational definition which stipulates when a rule can be applied by using sub-terms and substitution.  The second is an inductive definition which constructs the rewrite relation using generators and closure operations. When rewriting theory expanded to graphs in the 1970's, only the operational definition prevailed. The double pushout mechanism performed the substitution.  Then in the 1998, Gadducci and Heckel introduced an inductive definition to graph rewriting \cite{Gadd_IndGraphTrans}, thus allowing for analyses using structural induction. With the new technology of structured cospans, we can use their ideas to bring the inductive viewpoint to rewriting in a broader class of adhesive categories beyond directed graphs. 

A first key idea in developing the inductive definition of graph rewriting was to use open graphs. To \emph{openify} objects beyond graphs, we use structured cospans.

A second key idea is using an equivalence between two classes of rewrite relations.  In the context of graph rewriting, this result states that the rewrite relation for a graph grammar $(\Graph , \{ \ell_j \gets k_j \to r_j\})$ is the same as for its underlying \emph{discrete graph grammar} $( \Graph , \{ \ell_j \gets \flat k_j \to r_j\})$ where we replace the graphs $k_j$ with their underlying discrete graphs $\flat k_j$ and restrict the maps accordingly. To extend this to rewriting in adhesive categories, we start with an adjunction $L \dashv R$ with a monic counit and interpret the comonad $ \flat \bydef LR $ as sending an object $x$ to the \emph{discrete} object $ \flat x $ underlying $ x $.

In summation, we have the necessary ingredients to pursue an inductive viewpoint of rewriting: structured cospans provide objects with interfaces and the comonad $\flat$ provides discrete grammars. With these, we generalize Gadducci and Heckel's inductive-style rewriting on graphs to a broader collection of objects (Theorem \ref{thm:inductive-rewriting}), including hypergraphs and Petri Nets. 

\subsection{Outline and Contributions}
\label{sec:outl-contr}

Section \ref{sec:background} contains a brief overview of background material on double pushout rewriting and structured cospans.  

Section \ref{sec:str-csp-and-rewriting} introduces rewriting to structured cospans.  Section \ref{sec:characterizing-str-csp-as-adhesive-topoi} contains the main result of this paper; after introducing a new category $ _L\StrCsp $ with structured cospans as objects, we show that this category is adhesive (Theorem \ref{thm:strcsp-is-adhes}) and, under slightly stronger conditions, a topos (Theorem \ref{thm:strcsp-istopos}) constructed functorially in $ L $ (Theorem \ref{thm:strcsp-isfunctorial}). After discussing structured cospan grammars and their languages in Section \ref{sec:RewritingStrCsp}, we briefly mention the properties that they then inherit in Section \ref{sec:basic-properties}.

Section \ref{sec:an-application} uses structured cospans to bring an inductive viewpoint of rewriting to a broader class of objects. A brief history of the inductive viewpoint is provided in Section \ref{sec:hist-inductive-rewriting}. Discrete grammars are introduced in Section \ref{sec:discrete-grammars} followed, in Section \ref{sec:grammar-and-discrete}, by a sufficient condition for a grammar and its discrete grammar to induce the same language (Lemma \ref{thm:po-comp-same-rewrite-rel}). This condition is shown to hold for directed hypergraphs, undirected hypergraphs, Petri nets, marked Petri nets, and typed objects as constructed using slice categories.  Finally, the main application of this paper is given in Section \ref{sec:inductive-rewriting-hypergraphs}: the construction of the inductive viewpoint (Theorem \ref{thm:inductive-rewriting}).

\subsection{Related Work}
\label{sec:related-work}

The theory of rewriting is a mature field, particularly the rewriting of various kinds of graphs and the axiomatization of graph rewriting via adhesive categories.  

The concept of open networks has been previously explored, including open graphs, open hypergraphs and open Petri nets \cite{ baldan-etal-comp-modeling-open-nets, DixKiss_OpenGraphs, sassone-sobocinski_congruence-graph-rewriting}. Approaches included various techniques to equip graph-type objects with an interface.  Structured cospans provided a unified technique that is applicable to a broader class of open networks \cite{StructuredCospans, StructuredDecorated, baez-master_open-pnets, Cic_SpCsp}.  An alternative approach is Fong's decorated cospans \cite{DecorCsp}, which we do not address here. Rewriting open graphs, hypergraphs, Petri nets, and string diagrams have been studied \cite{bonchi-kissinger_rewriting-mod-sms,bonchi-etal_graph-rewriting-interfaces,ehrig_deriving-bisimulation, sassone-sobocinski_congruence-graph-rewriting, sassone-sobocinski_congruence-pnets} but this paper is the first to introduce rewriting specifically to structured cospans. This contribution places previous efforts into a unified framework, using general categorical techniques.

In addition to the above, this paper draws on the theory of adhesive categories \cite{LackSobo_Adhesive}, and the presheaf perspective on both hypergraphs \cite{bonchi-etal_graph-rewriting-interfaces} and Petri nets \cite{johnstone_quasi}. 


\section{Background Material: DPO Rewriting and Structured Cospans}
\label{sec:background}

This section contains background material on our main topics: double pushout rewriting and structured cospans.  Nothing novel is presented here.

\subsection{Double Pushout Rewriting}
\label{sec:dpo-rewriting}

Double pushout rewriting has an established literature, so we use this section to cover the fundamentals and to establish our conventions. The interested reader can see Ehrig, et.~al.~ \cite{Ehrig_GraphGram} to learn about graph rewriting or Lack and Soboci\'{n}ski \cite{LackSobo_Adhesive} for an axiomatic approach based on adhesive categories. 

Rewriting starts with the notion of a \df{rewrite rule}, or simply \df{rule}. This is a span $\ell \gets k \to r$ with two monic arrows. The interpretation of this rule is that $ \ell $ can be replaced by $ r $ and $ k $ is the part of $ \ell $ that does not change. A \df{grammar} $(\C,P)$ is defined to be an adhesive category $ \C $ together with a finite set of rules  $ P \bydef \{ \spn{\ell_j}{k_j}{r_j} \} $. 

We can apply a rule $ \ell \gets k \to r $ from $P$ to an object $ \ell' $ of $\C$ using any arrow $ m \from \ell \to \ell' $ for which there exists a pushout complement, that is an object $ k' $ fitting into a pushout diagram
\[
  \begin{tikzpicture}
    \node (l) at (0,1) {$ \ell $};
    \node (k) at (1,1) {$ k $};
    \node (g) at (0,0) {$ \ell' $};
    \node (d) at (1,0) {$ k' $};
    \draw [cd]
    (k) edge node[]{$  $}      (l)
    (k) edge node[]{$  $}      (d)
    (l) edge node[left]{$ m $} (g)
    (d) edge node[]{$  $}      (g); 
    \draw (0.3,0.4) -- (0.4,0.4) -- (0.4,0.3);
  \end{tikzpicture}
\]
A pushout complement need not exist, but if it does and the map $ k \to \ell $ is monic, then it is unique up to isomorphism \cite[Lem.~15]{LackSobo_Adhesive}.

Every application of a rule begets a new rule. Applying $ \ell \gets k \to r $ to $\ell'$ along $ m \from \ell \to \ell' $ induces a \df{derived rule} $ \ell' \gets k' \to r' $ obtained as the bottom row of the double pushout diagram
\[
  \begin{tikzpicture}
    \node (1) at (0,1) {$ \ell $};
    \node (2) at (1,1) {$ k $};
    \node (3) at (2,1) {$ r $};
    \node (4) at (0,0) {$ \ell' $};
    \node (5) at (1,0) {$ k' $};
    \node (6) at (2,0) {$ r' $};
    \draw [cd]
    (2) edge (1)
    (2) edge (3)
    (5) edge (4)
    (5) edge (6)
    (1) edge node[left]{$ m $} (4)
    (2) edge (5)
    (3) edge (6); 
    \draw (0.3,0.4) -- (0.4,0.4) -- (0.4,0.3);
    \draw (1.7,0.4) -- (1.6,0.4) -- (1.6,0.3);
  \end{tikzpicture}
\]
This diagram expresses a three-stage process whereby $ m $ selects a copy of $ \ell $ inside $ \ell' $, this copy is replaced by $ r $, and the resulting object $ r' $ is returned.  Because pushouts preserve monos in adhesive categories, a derived rule is, in fact, a rule.

A grammar $( \C,P )$ induces a collection $dP$ of all derived rules obtained by applying a rule in $P$ to an object in $\C$.  We can use $dP$ to analyze the grammar $(\C,P)$ by constructing the \emph{rewrite relation} $\deriv{}{}$. The meaning of $\deriv{x}{y}$ is that we can rewrite $x$ into $y$ by applying a sequence of rules in $dP$. To precisely define the rewrite relation, we start by constructing a relation $ \dderiv{}{} $ on the objects of $ \C $ by setting $ \dderiv{\ell'}{r'} $ if there exists a rule $ \spn{\ell'}{k'}{r'} $ in $ dP $. However, $ \dderiv{}{} $ does not capture enough information about $( \C,P )$, which is why we define the \df{rewrite relation} $ \deriv{}{} $ to be the reflexive and transitive closure of $ \dderiv{}{} $.

The rewrite relation can be encoded as arrows in a category, which we call the \df{language of a grammar} $ \Lang (\C,P) $, whose objects are those of $\C$ and arrows are generated by $\dderiv{}{}$. The language is constructed so that there is an arrow $ x \to y $ in $ \Lang (\C,P) $ if and only if $ \deriv{x}{y} $. 

In the rewriting literature, the terms ``language'' and ``rewrite relation'' are often interchangeable. However, we give them slightly different meanings in order to help orient the reader. Namely, we use ``rewrite relation'' when giving a relational perspective and ``language'' when giving the category theoretical perspective. 


Though there is more to the theory of rewriting than is provided in
this section, we have developed enough of the theory to continue our goal of introducing rewriting to structured cospans.

\subsection{Structured Cospans}
\label{sec:structured-cospans}

Baez and Courser \cite{StructuredCospans} introduced structured cospans as a framework to study open networks.  A network is \df{open} when equipped with a mechanism by which it can connect to any compatible network.  For example, a vacuum cleaner can connect with the electrical grid via an electrical socket. A pulley network can connect to a mechanical motor. An open network stands in contrast to a closed network that cannot interact with its outside environment.

\begin{definition}[Structured Cospans]
  \label{df:str-csp}
  Given a functor $L \from \A \to \X$ where $\X$ has pushouts, an \df{$L$-structured cospan}, or simply \df{structured cospan}, is a diagram $\csp{La}{x}{Lb}$ in $\X$.  
\end{definition}

To interpret $ \csp{La}{x}{Lb} $ as an open network, take $x$ to represent the network with inputs $ La $ chosen by the arrow $ La \to x$ and outputs $ Lb $ chosen by $ x \gets Lb $. This open network can now connect to any other open network with inputs $ Lb $, say $ \csp{Lb}{y}{Lc} $. We form the \emph{composite} of the two open networks by connecting $ x $ to $ y $ along their common interface $ Lb $. Mathematically, this amounts to taking the pushout of $ x $ and $ y $ over $ Lb $, thus giving the composite network $ \csp{La}{x+_{Lb}y}{Lc} $. We capture this in a category $_{L} \Csp$ whose objects are those of $\A$ and arrows $a \to b$ are structured cospans $\csp{La}{x}{Lb}$ with composition given by pushout. 

\begin{example} \label{ex:open_graphs}
  Structured cospans can be used to slightly generalize open graphs.  Set theoretically, a graph is open when equipped with two subsets of its nodes, one set serving as inputs and the other as outputs. When the inputs of one open graph coincide with the outputs of another, they can be composed. For example, the pair of open graphs
  
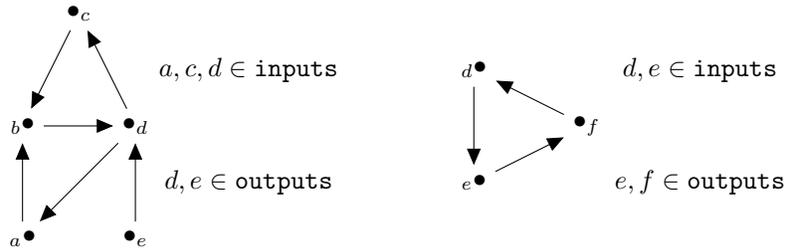
\[
    \begin{tikzpicture}[scale=0.6]
      \begin{scope}
        \node (1) at (0,0) {{$ _a \bullet $}};
        \node (2) at (0,2) {{$ _b \bullet $}};
        \node (3) at (2,2) {{$ \bullet_d $}};
        \node (4) at (2,0) {{$ \bullet_e $}};
        \node (5) at (1,4) {{$ \bullet_c $}};
        \draw [graph]
        (1) edge[] (2)
        (2) edge[] (3)
        (3) edge[] (1)
        (4) edge[] (3)
        (3) edge[] (5)
        (5) edge[] (2);
        \node () at (-2.5,3) {$ a,c,d \in \mathtt{ inputs} $ };
        \node () at (-2.5,1) {$ d,e \in \mathtt{ outputs} $};
      \end{scope}
      \begin{scope}[shift={(5,0)}]
        \node (1) at (0,3) {{$ _{d} \bullet $}};
        \node (2) at (0,1) {{$ _{e} \bullet $}};
        \node (3) at (2,2) {{$ \bullet_{f} $}};
        \draw [graph] 
        (1) edge[] (2)
        (2) edge[] (3)
        (3) edge[] (1); 
        \node () at (4.5,3) {$ d,e \in \mathtt{ inputs} $};
        \node () at (4.5,1) {$ e,f \in \mathtt{ outputs} $};  
      \end{scope}
    \end{tikzpicture}
  \]
  compose by gluing the corresponding nodes together, forming the new open graph
  
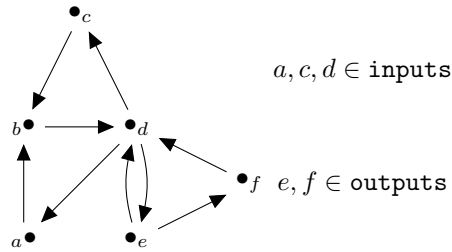
\[
    \begin{tikzpicture}[scale=0.6]
      \node (1) at (0,0) {{$ _a \bullet $}};
      \node (2) at (0,2) {{$ _b \bullet $}};
      \node (3) at (2,2) {{$ \bullet_{d} $}};
      \node (4) at (2,0) {{$ \bullet_{e} $}};
      \node (5) at (1,4) {{$ \bullet_c $}};
      \node (6) at (4,1) {{$ \bullet_f$}};
      \draw [graph]
      (1) edge[] (2)
      (2) edge[] (3)
      (3) edge[] (1)
      (3) edge[bend left=15] (4)
      (3) edge[] (5)
      (5) edge[] (2)
      (4) edge[bend left=15] (3)
      (4) edge[] (6)
      (6) edge[] (3); 
      \node () at (6.5,3) {$ a,c,d \in \type{inputs} $};
      \node () at (6.5,1) {$ e,f \in \type{outputs}$};
    \end{tikzpicture}  
  \]
  To define an open graph as a structured cospan, consider $$\adjunction{\Graph}{L}{R}{\Set}$$ where $ L $ is the discrete graph functor and $ R $ forgets the graph edges. The above open graphs can be presented as the structured cospans
  
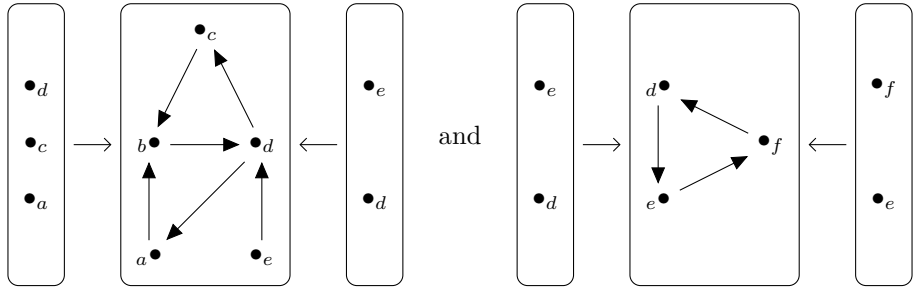
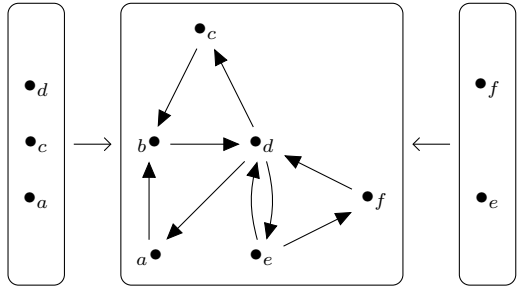
\[
    \begin{tikzpicture}[
      baseline=(current bounding box.center),
      scale=0.6]
      \begin{scope}
        \node () at (0,1) {$ \bullet_a $};
        \node () at (0,2) {$ \bullet_c $};
        \node () at (0,3) {$ \bullet_d $};
        \draw [rounded corners]
        (-0.5,-0.5) rectangle (0.5,4.5);
        \node (InDom) at (0.5,2) {};
      \end{scope}
      \begin{scope}[shift={(2,0)}]
        \node (1) at (0,0) {{$ _a \bullet $}};
        \node (2) at (0,2) {{$ _b \bullet $}};
        \node (3) at (2,2) {{$ \bullet_d $}};
        \node (4) at (2,0) {{$ \bullet_e $}};
        \node (5) at (1,4) {{$ \bullet_c $}};
        \draw [graph]
        (1) edge[] (2)
        (2) edge[] (3)
        (3) edge[] (1)
        (4) edge[] (3)
        (3) edge[] (5)
        (5) edge[] (2);
        \draw [rounded corners]
        (-0.5,-0.5) rectangle (2.5,4.5);
        \node (InCod) at (-0.5,2) {};
        \node (OutCod) at (2.5,2) {};
      \end{scope}
      \begin{scope}[shift={(6,0)}]
        \node () at (0,1) {$ \bullet_e $};
        \node () at (0,3) {$ \bullet_d $};
        \draw [rounded corners]
        (-0.5,-0.5) rectangle (0.5,4.5);
        \node (OutDom) at (-0.5,2) {};
      \end{scope}
      \draw [cd]
      (InDom) edge[] (InCod)
      (OutDom) edge[] (OutCod);
    \end{tikzpicture}
    \quad \text{ and } \quad
    \begin{tikzpicture}[
      baseline=(current bounding box.center),
      scale=0.6]
      \begin{scope}
        \node () at (0,1) {$ \bullet_e $};
        \node () at (0,3) {$ \bullet_d $};
        \draw [rounded corners]
        (-0.5,-0.5) rectangle (0.5,4.5);
        \node (InDom) at (0.5,2) {};
      \end{scope}
      \begin{scope}[shift={(2,0)}]
        \node (1) at (0,3) {{$ _{d} \bullet $}};
        \node (2) at (0,1) {{$ _{e} \bullet $}};
        \node (3) at (2,2) {{$ \bullet_{f} $}};
        \draw [graph] 
        (1) edge[] (2)
        (2) edge[] (3)
        (3) edge[] (1); 
        \draw [rounded corners]
        (-0.5,-0.5) rectangle (2.5,4.5);
        \node (InCod) at (-0.5,2) {};
        \node (OutCod) at (2.5,2) {};
      \end{scope}
      \begin{scope}[shift={(6,0)}]
        \node () at (0,1) {$ \bullet_e $};
        \node () at (0,3) {$ \bullet_f $};
        \draw [rounded corners]
        (-0.5,-0.5) rectangle (0.5,4.5);
        \node (OutDom) at (-0.5,2) {};
      \end{scope}
      \draw [cd]
      (InDom) edge[] (InCod)
      (OutDom) edge[] (OutCod);
    \end{tikzpicture}
  \]
  with composite
  \[
    \begin{tikzpicture}[scale=0.6]
      \begin{scope}
        \node () at (0,1) {$ \bullet_a $};
        \node () at (0,2) {$ \bullet_c $};
        \node () at (0,3) {$ \bullet_d $};
        \draw [rounded corners]
        (-0.5,-0.5) rectangle (0.5,4.5);
        \node (InDom) at (0.5,2) {};
      \end{scope}
      \begin{scope}[shift={(2,0)}]
        \node (1) at (0,0) {{$ _a \bullet $}};
        \node (2) at (0,2) {{$ _b \bullet $}};
        \node (3) at (2,2) {{$ \bullet_{d} $}};
        \node (4) at (2,0) {{$ \bullet_{e} $}};
        \node (5) at (1,4) {{$ \bullet_c $}};
        \node (6) at (4,1) {{$ \bullet_f$}};
        \draw [graph]
        (1) edge[] (2)
        (2) edge[] (3)
        (3) edge[] (1)
        (3) edge[bend left=15] (4)
        (3) edge[] (5)
        (5) edge[] (2)
        (4) edge[bend left=15] (3)
        (4) edge[] (6)
        (6) edge[] (3);  
        \draw [rounded corners]
        (-0.5,-0.5) rectangle (4.5,4.5);
        \node (InCod) at (-0.5,2) {};
        \node (OutCod) at (4.5,2) {};
      \end{scope}
      \begin{scope}[shift={(8,0)}]
        \node () at (0,1) {$ \bullet_e $};
        \node () at (0,3) {$ \bullet_f $};
        \draw [rounded corners]
        (-0.5,-0.5) rectangle (0.5,4.5);
        \node (OutDom) at (-0.5,2) {};
      \end{scope}
      \draw [cd]
      (InDom) edge[] (InCod)
      (OutDom) edge[] (OutCod);
    \end{tikzpicture}
  \]
  Note, that the functions forming the legs of the cospan need not be monic, hence calling the structured cospan version of open graphs a slight generalization of the set theoretic version of open graphs. In this example, the category $ _L\Csp $ has sets for objects and open graphs for arrows.
\end{example}

\section{Structured Cospans and Rewriting}
\label{sec:str-csp-and-rewriting}

In this section, we give our main result, that structured cospans form an adhesive category or a topos under mild conditions (Section \ref{sec:characterizing-str-csp-as-adhesive-topoi}). It follows from this result that structured cospans admit a rewriting theory.  We then define a grammar and its language for structured cospans (Section \ref{sec:RewritingStrCsp}). These require a bit more discussion than simply importing the ideas from the current rewriting literature. Finally, we reflect on some of the properties that rewriting structured cospans will inherit, simply as a matter of being adhesive (Section \ref{sec:basic-properties}).

\subsection{Characterizing Structured Cospan Categories as Adhesive or Topoi}
\label{sec:characterizing-str-csp-as-adhesive-topoi}

Every adhesive category supports a rich rewriting theory. This fact underpins our efforts to introduce rewriting to structured cospans. And so, in this section, we reintroduce the notion of structured cospans in a different context and find sufficient conditions for them to be adhesive.

Recall that, in Section \ref{sec:background}, we saw that rewriting operates on the objects of an adhesive category, not the arrows. Therefore, we cannot hope to rewrite structured cospans inside the category $ _L\Csp $. Our task, now, is to build a category where structured cospans are objects and then to show that it is adhesive. 

\begin{definition} \label{df:morph-of-strcsp}
  Fix an adjunction $ L \dashv R \from \X \to \A $. Define $ _L \StrCsp $ to be the category whose objects are structured cospans and arrows from $ La \to x \gets Lb $ to $ La' \to x' \gets Lb' $ are triples of arrows $ ( f,g,h ) $ fitting into the commuting diagram
  \[
    \begin{tikzpicture}
      \node (1) at (0,1) {\( La \)};
      \node (2) at (1,1) {\( x \)};
      \node (3) at (2,1) {\( Lb \)};
      \node (4) at (0,0) {\( La'\)};
      \node (5) at (1,0) {\( x' \)};
      \node (6) at (2,0) {\( Lb' \)};
      \draw [cd]
      (1) edge node[]{} (2)
      (3) edge node[]{} (2)
      (4) edge node[]{} (5)
      (6) edge node[]{} (5)
      (1) edge node[left]{$Lf$}  (4)
      (2) edge node[left]{$g$}   (5)
      (3) edge node[right]{$Lh$} (6);
    \end{tikzpicture}
  \]
\end{definition}

Now, in order to enable the rewriting of structured cospans in the sense of Section {\ref{sec:background}}, $_{L} \StrCsp$ must be adhesive.  To achieve this, we outline two sets of conditions in order to later prove results of different strengths.

\begin{definition} \label{df:conditions}
  Fix an adjunction $ L \dashv R \from \X \to \A $. We say that this adjunction satisfies the \df{adhesivity condition} if $\X$ and $\A$ are adhesive and $L$ preserves pullbacks.  We say that this adjunctions satisfies the \df{topos condition} if $R$ is a geometric morphism between topoi.  
\end{definition}

Every topos is adhesive \cite{LackSobo_TopsIsAdh}, meaning the topos condition is strictly stronger than the adhesive condition. As shown in the next two results, when structured cospans are built with $L \dashv R$ satisfying the adhesivity condition, then $_{L} \StrCsp$ is adhesive.  Similarly, $_{L} \StrCsp$ is a topos when $L \dashv R$ satisfies the topos condition.

\begin{therm} \label{thm:strcsp-is-adhes}
  If $ L \dashv R \from \X \to \A $ satisfies the adhesive condition, then $ _{L}\StrCsp $ is adhesive.
\end{therm}

\begin{proof}
  Pullbacks exist in $_{L} \StrCsp$ because both $\X$ and $\A$ have pullbacks and $L$ preserves them.  The same is true about pushouts of monos.  That a pushout of a monic map is a Van Kampen square follows from this property holding in both $\X$ and $\A$ and $L$ preserving pullbacks and pushouts.
\end{proof}

While in this paper, we will largely work with the generality of adhesive categories, many objects of rewriting interest actually form topoi.  Graphs are the primary example. Therefore, we show that under slightly stronger conditions, structured cospans form a topos and that this construction is functorial in $L$.

\begin{therm} \label{thm:strcsp-istopos}
  Let $ L \dashv R \from \X \to \A $ be a geometric morphism.  The category $ _{L}\StrCsp $ is a topos.
\end{therm}

\begin{proof}
  By adjointness, $ _{L}\StrCsp $ is equivalent to the category whose objects are cospans of form $ a \to Rx \gets b $ and morphisms are triples $ ( f,g,h ) $ fitting into the commuting diagram
  \[
    \begin{tikzpicture}
      \node (1) at (0,1) {\( a \)};
      \node (2) at (1,1) {\( Rx \)};
      \node (3) at (2,1) {\( b \)};
      \node (4) at (0,0) {\( a' \)};
      \node (5) at (1,0) {\( Rx' \)};
      \node (6) at (2,0) {\( b' \)};
      \draw [cd] (1) to  node [] {} (2);
      \draw [cd] (3) to node [] {} (2);
      \draw [cd] (4) to node [] {} (5);
      \draw [cd] (6) to node [] {} (5);
      \draw [cd] (1) to node [left] {$f$} (4);
      \draw [cd] (2) to node [left] {$Rg$} (5);
      \draw [cd] (3) to node [left] {$h$} (6); 
    \end{tikzpicture}
  \]
  This, in turn, is equivalent to the comma category $ ( \A \times \A ) \downarrow \Delta R  $ where $ \Delta \from \A \to \A \times \A $ is the diagonal functor, a right adjoint. Because $ \Delta R $ is a right adjoint, $ ( \A \times \A ) \downarrow \Delta R $ is an Artin gluing \cite{Wraith_ArtinGlue}, therefore a topos.
\end{proof}

\begin{therm} \label{thm:strcsp-isfunctorial}
  Denote by $ \Topos $ the category of topoi and geometric morphisms. There is a contravariant functor $ _{(-)}\StrCsp \from [\bullet \to \bullet, \Topos ] \to  \Topos$ defined by
  \[
    \begin{tikzpicture}
      \begin{scope}
        \node (1) at (-1,1) {$\X$};
        \node (2) at (-1,-1) {$\X'$};
        \node (3) at (1,1) {$\A$};
        \node (4) at (1,-1) {$\A'$};
        \draw [transform canvas = {xshift=1ex}, ->] (1) edge node [right] {$G$} (2);
        \draw [transform canvas = {xshift=-1ex}, ->] (2)  edge node [left]  {$F$}  (1);
        \draw [transform canvas = {yshift=1ex}, ->] (3)  edge node [above] {$L$}  (1);
        \draw [transform canvas = {yshift=-1ex}, ->] (1)  edge node [below] {$R$}  (3);
        \draw [transform canvas = {yshift=-1ex}, ->] (2)   edge node [below] {$R'$} (4);
        \draw [transform canvas = {yshift=1ex}, ->] (4) edge node [above] {$L'$} (2);
        \draw [transform canvas = {xshift=-1ex}, ->] (4)   edge node [left] {$F'$} (3);
        \draw [transform canvas = {xshift=1ex}, ->] (3) edge node [right]  {$G'$} (4);
        \node (5) at (0,-1) {\scriptsize{\( \perp \)}};
        \node (6) at (0,1) {\scriptsize{\( \perp \)}};
        \node (7) at (-1,0) {\scriptsize{\( \dashv \)}};
        \node (8) at (1,0) {\scriptsize{\( \dashv \)}};
      \end{scope}
      \begin{scope}[shift={(3,0)}]
        \node (1) at (0,0) { $\xmapsto{ _{(-)}\StrCsp }$ };
      \end{scope}
      \begin{scope}[shift={(5,0)}]
        \node (1) [] at (0,0) {\( _{L'}\StrCsp \)};
        \node (2) [] at (2,0) {\( _{L}\StrCsp \)};
        \node (3) at (1,0) {\scriptsize{ \( \perp \) }};
        \draw [cd]
        ([yshift=1pt]2.west) to
        node [above] {$\Theta$}
        ([yshift= 1pt]1.east);
        \draw [cd]
        ([yshift= -1pt]1.east) to
        node [below] {$\Psi$}
        ([yshift= -1pt]2.west);  
      \end{scope}
    \end{tikzpicture}
  \]
  which is in turn given by
  \[
    \begin{tikzpicture}
      \begin{scope}
        \node (1) at (0,1) {\( La \)};
        \node (2) at (1,1) {\( x \)};
        \node (3) at (2,1) {\( Lb \)};
        \node (4) at (0,0) {\( Lc \)};
        \node (5) at (1,0) {\( y \)};
        \node (6) at (2,0) {\( Ld \)};
        \draw [cd]
        (1) edge node [above] {$ m $}  (2)
        (3) edge node [above] {$ n $}  (2)
        (4) edge node [below] {$ o $}  (5)
        (6) edge node [below] {$ p $}  (5)
        (1) edge node [left]  {$ Lf $} (4)
        (2) edge node [left]  {$ g $}  (5)
        (3) edge node [right]  {$ Lh $} (6);
      \end{scope}
      \begin{scope}[shift={(4,0)}]
        \node (1) at (0,0.5) { $ \xmapsto{ \Theta } $ };
      \end{scope}
      \begin{scope}[shift={(6,0)}]
        \node (1) at (0,1) {\( L'G'a \)};
        \node (2) at (1.5,1) {\( Gx \)};
        \node (3) at (3,1) {\( L'G'b \)};
        \node (4) at (0,0) {\( L'G'c \)};
        \node (5) at (1.5,0) {\( Gy \)};
        \node (6) at (3,0) {\( L'G'd \)};
        \draw [->] (1) to node [above] {\scriptsize{\( Gm \)}} (2);
        \draw [->] (3) to node [above] {\scriptsize{\( Gn \)}} (2);
        \draw [->] (4) to node [below] {\scriptsize{\( Go \)}} (5);
        \draw [->] (6) to node [below] {\scriptsize{\( Gp \)}} (5);
        \draw [->] (1) to node [left] {\scriptsize{\( L'G'f \)}} (4);
        \draw [->] (2) to node [left] {\scriptsize{\( Gg \)}} (5);
        \draw [->] (3) to node [right] {\scriptsize{\( L'G'h \)}} (6);  
      \end{scope}
    \end{tikzpicture}
  \]
  \[
    \begin{tikzpicture}
      \begin{scope}
        \node (1) at (0,1) {\( L'a' \)};
        \node (2) at (1,1) {\( x' \)};
        \node (3) at (2,1) {\( L'b' \)};
        \node (4) at (0,0) {\( L'c' \)};
        \node (5) at (1,0) {\( y' \)};
        \node (6) at (2,0) {\( L'd' \)};
        \draw [->] (1) edge node [above] {\scriptsize{\( m' \)}} (2);
        \draw [->] (3) edge  node [above] {\scriptsize{\( n' \)}} (2);
        \draw [->] (4) edge node [below] {\scriptsize{\( o' \)}} (5);
        \draw [->] (6) edge node [below] {\scriptsize{\( p' \)}} (5);
        \draw [->] (1) edge node [left] {\scriptsize{\( L'f' \)}} (4);
        \draw [->] (2) edge node [left] {\scriptsize{\( g' \)}} (5);
        \draw [->] (3) edge node [right] {\scriptsize{\( L'h' \)}} (6);
      \end{scope}
      \begin{scope}[shift={(4,0)}]
        \node (1) at (0,0.5) { $ \xmapsto{ \Psi } $ };
      \end{scope}
      \begin{scope}[shift={(6,0)}]
        \node (1) at (0,1) {\( LF'a' \)};
        \node (2) at (1.5,1) {\( Fx' \)};
        \node (3) at (3,1) {\( LF'b' \)};
        \node (4) at (0,0) {\( LF'c' \)};
        \node (5) at (1.5,0) {\( Fy' \)};
        \node (6) at (3,0) {\( LF'd' \)};
        \draw [->] (1) edge node [above] {\scriptsize{\( Fm' \)}} (2);
        \draw [->] (3) edge node [above] {\scriptsize{\( Fn' \)}} (2);
        \draw [->] (4) edge node [below] {\scriptsize{\( Fo' \)}} (5);
        \draw [->] (6) edge node [below] {\scriptsize{\( Fp' \)}} (5);
        \draw [->] (1) edge node [left] {\scriptsize{\( LF'f' \)}} (4);
        \draw [->] (2) edge node [left] {\scriptsize{\( Fg' \)}} (5);
        \draw [->] (3) edge node [right] {\scriptsize{\( LF'h' \)}} (6);  
      \end{scope}
    \end{tikzpicture}
  \]
\end{therm}

\begin{proof}
  In light of Theorem \ref{thm:strcsp-istopos}, it suffices to show that $ \Theta \dashv \Psi $ gives a geometric  morphism. Let $ \ell $ and $\ell'$ respectively denote the $ L $-structured cospan $ La \to x \gets Lb$ and the $ L' $-structured cospan $ L'a' \to x' \gets L'b'$. Denote the unit and counit for $F \dashv G$ by $ \eta $, $ \varepsilon $ and for $ F' \dashv G' $ by $ \eta' $, $ \varepsilon' $.  The assignments
  \begin{align*}
    \left( (f,g,h)\from\ell\to\Psi\ell' \right)
      & \mapsto
        \left( (\epsilon'\circ G'f,
        \epsilon\circ Gg,
        \epsilon'\circ G'h )
        \from\Theta\ell\to\ell'\right) \\
    \left( (f',g',h')\from\Theta\ell\to\ell'\right)
      & \mapsto
        \left( (Ff'\circ\eta',F'g'\circ\eta,F'h'\circ\eta')
        \from\ell\to\Theta'\ell'\right) 
  \end{align*}  
  give a bijection $ \hom (\Theta\ell,\ell') \simeq \hom (\ell,\Psi\ell' )$. The naturality in $\ell$ and $\ell'$ follow from the naturality of $\eta$, $\epsilon$, $\eta'$, and $\epsilon'$. The left adjoint $\Theta$ preserves finite limits because they are taken pointwise and $ L' $, $ G $, and $ G' $ all preserve finite limits. Identity and composition are easily checked to be preserved.  
\end{proof}

For the remainder of the paper, we work with structured cospans that satisfy the adhesive condition unless otherwise stated.

\subsection{Rewriting structured cospans}
\label{sec:RewritingStrCsp}

We now know that the category $ _L\StrCsp $ of structured cospans and their morphisms are adhesive and, therefore, supports a rich rewriting theory. In this section, we begin to develop this theory.

\begin{definition} \label{df:rewrite-rule-str-csp}
  A \df{rewrite rule of structured cospans} is an isomorphism class of spans of structured cospans \cite{Cic_SpCsp, CicCour_SpCspTopos} of the form
  \[
    \begin{tikzpicture}[baseline=(current bounding box.center)]
      \node (x') at (1,2) {$ x' $};
      \node (La) at (0,1) {$ La $};
      \node (x) at (1,1) {$ x $};
      \node (Lb) at (2,1) {$ Lb $};
      \node (x'') at (1,0) {$ x'' $};
      \draw [>->] (x) to (x');
      \draw [>->] (x) to (x'');
      \draw[cd]
      (La) edge[] (x)
      (Lb) edge[] (x)
      (La) edge[] (x')
      (Lb) edge[] (x')
      (La) edge[] (x'')
      (Lb) edge[] (x'');
    \end{tikzpicture}
  \]
  The marked arrows $ \rightarrowtail $ are monic.
\end{definition}

The conceit of this rule is that the structured cospan $ \csp{La}{x'}{Lb} $ along the top of the diagram is rewritten as the structured cospan $ \csp{La}{x''}{Lb} $ along the bottom. Here, we are orienting our diagrams as `top rewrites to bottom' instead of `left rewrites to right'.  A non-superficial difference is that this diagram is not merely a span in $ _L\StrCsp $ with monic legs.  We disallow the inputs and outputs from being rewritten.

In order to apply this rule to other graphs, we need the notions of grammars and derived rules for structured cospans.

\begin{definition}
  A \df{grammar of structured cospans} is a pair $(_L\StrCsp, P)$ where $P$ is a set of rewrite rules of structured cospans.
\end{definition}

\begin{proposition}
  Derived rules of structured cospan rewrites are structured cospan rewrite rules.
\end{proposition}

\begin{proof}
  Consider the diagram
  \[
    \begin{tikzpicture}[scale=0.75, baseline=(current  bounding  box.center)]
      \node (x')  at (-2,4) {$ x' $};
      \node (y')  at (-2,1) {$ y' $};
      \node (La) at (0,3) {$ La $};
      \node (x)  at (1,4) {$ x $};
      \node (Lb) at (2,5) {$ Lb $};
      \node (Lc) at (0,0) {$ Lc $};
      \node (y)  at (1,1) {$ y $};
      \node (Ld) at (2,2) {$ Ld $};
      \node (x'')  at (4,4) {$ x'' $};
      \node (y'')  at (4,1) {$ y'' $};
      \draw [cd]
      (Lb)   edge[] (x')
      (Lb)   edge[] (x'') 
      (Lb)   edge[] (Ld)
      (Ld)   edge[] (y')
      (Ld)   edge[] (y''); 
      \draw [white,line width=2mm]
      (x) -- (x')
      (x) -- (x'') 
      (x') -- (y')
      (x) -- (y)
      (x'') -- (y'')
      (y) -- (y')
      (y) -- (y'');
      \draw [cd]
      (x)   edge[>->] (x')
      (x)   edge[>->] (x'') 
      (x')  edge[] (y')
      (x)   edge[] (y)
      (x'') edge[] (y'')
      (y)   edge[>->] (y')
      (y)   edge[>->] (y''); 
      \draw [white,line width=1mm]
      (La) -- (x')
      (La) -- (x'') 
      (La)   -- (Lc)
      (Lc)   -- (y')
      (Lc)   -- (y'');
      \draw [cd]
      (La)  edge[] (x')
      (La)  edge[] (x'')      
      (La)  edge[] (Lc)
      (Lc)   edge[] (y')
      (Lc)   edge[] (y'');  
      \draw [cd]
      (Lb)   edge[] (x)
      (La)   edge[] (x)
      (Ld)   edge[] (y)
      (Lc)   edge[] (y);
    \end{tikzpicture}
  \]
  where the top face is a rewrite rule of structured cospans, and the back left and front left faces give a matching map, with pushout complement, of structured cospans. Then the bottom face is also a rewrite rule of structured cospans because pushouts preserve monos in an adhesive category.   
\end{proof}

Now that we are certain that derivation preserves rewrite rules of structured cospans, we can construct the language. Similar to the construction in Gadducci and Heckel \cite{Gadd_IndGraphTrans}, we realize the language as the free 2-category on a computad whose 0-cells are interface types, 1-cells are structured cospans, and 2-cells are generated by the rewrite rules of a given grammar.

\begin{definition}\label{df:lang}
  Given an adjunction $L \dashv R \from X \to A$ satisfying the adhesivity condition (Def.~\ref{df:conditions}) and a structured cospan grammar $(_L\StrCsp, P)$, define a computad $\Comp (_L\StrCsp, P)$ with $0$-cells the $\A$-objects, $1$-cells the arrows of $_L\StrCsp$ and, for each
  \[
    \begin{tikzpicture}[baseline=(current bounding box.center)]
      \node (x') at (1,2) {$ x' $};
      \node (La) at (0,1) {$ La $};
      \node (x) at (1,1) {$ x $};
      \node (Lb) at (2,1) {$ Lb $};
      \node (x'') at (1,0) {$ x'' $};
      \draw [>->] (x) to (x');
      \draw [>->] (x) to (x'');
      \draw[cd]
      (La) edge[] (x)
      (Lb) edge[] (x)
      (La) edge[] (x')
      (Lb) edge[] (x')
      (La) edge[] (x'')
      (Lb) edge[] (x'');
    \end{tikzpicture}
  \]
  derived from $P$, a $2$-cell $\gamma_p$
  \[
    \begin{tikzpicture}[
      baseline=(current bounding box.center), scale=0.75]
      \node (x') at (1,2) {$ x' $};
      \node (La) at (0,1) {$ La $};
      \node (Lb) at (2,1) {$ Lb $};
      \node (x'') at (1,0) {$ x'' $};
      \draw[cd]    
      (La) edge[] (x')
      (Lb) edge[] (x')
      (La) edge[] (x'')
      (Lb) edge[] (x'');
      \draw[double,->] (1,1.25) to (1,0.75);
    \end{tikzpicture}
  \]
  Apply to this computad the free 2-category functor, that is the left adjoint to the forgetful functor taking each 2-category to its underlying computad. The resulting 2-category is what we call $\LLang(_L\StrCsp, P)$.
\end{definition}

Note that we have two notions of language.  The first, written in non-bold type, $\Lang (\X,P)$ is a 1-category generated by the rewrite relation of a typical grammar. The second, written in bold type, $\LLang (_L\StrCsp, Q)$ is the 2-category generated by a computad for a structured cospan grammar. 

\subsection{Basic Properties Inherited}
\label{sec:basic-properties}

Thus far, we have constructed an adhesive category $_{L} \StrCsp$. Therefore, it enjoys the basic properties of adhesive categories, as detailed by Lack and Soboci\'{n}ski \cite{LackSobo_Adhesive}. Moreover, structured cospan grammars involve rewrite rules that are stronger than those for general adhesive grammars. For general adhesive categories, rewrite rules can be any span with monic legs while structured cospan rewrite rules are \emph{certain} spans with monic legs. Therefore, a structured cospan grammar is just a grammar in an adhesive category. It follows that the Local Church Rosser theorem and the Concurrency theorem hold for structured cospan grammars.  

\section{An Application: An Inductive Perspective of Rewriting}
\label{sec:an-application}

In this section, we shift our perspective from rewriting structured cospans to using structured cospans to provide an inductive viewpoint of rewriting. We start this section with a brief history of inductive rewriting, including its introduction to directed graphs. Then after establishing definitions, we close this section with our main result: sufficient conditions for categories to admit inductive rewriting on their objects. We illustrate this with directed and undirected hypergraphs.


\subsection{A Brief History of Inductive Rewriting}
\label{sec:hist-inductive-rewriting}

Before graph rewriting, there was both formal language rewriting and term rewriting. In these cases, there are two ways to define the rewrite relation.  The first way is called the \emph{operational method}, which applies a rule by substituting a sub-term for another term.  The second way is called the \emph{inductive method}, which constructs the rewrite relation using generators and closure operations. In classical graph rewriting, only the operational method existed, where substitution was achieved with the double pushout method.  Eventually, Gadducci and Heckel introduced an inductive method to construct the rewrite relation, opening the way to analyze graph grammars through structural induction. 

These authors accomplished this by defining what are now more commonly called \emph{open graphs} defined using a cospan $d_1 \to g \gets d_2$ where the images of discrete graphs $d_1$ and $d_2$ identified the inputs and outputs of a graph $g$.  This idea is a precursor for structured cospans. The inductively defined rewrite relation for a graph grammar used this open graph construction to encode the relation inside the hom-set of a particular 2-category. We emphasize here that the goal was not about rewriting open graphs, but to use open graphs to talk about rewriting graphs.

In this section, we adapt their ideas to give an inductive definition of the rewrite relation for a grammar $ ( \X,P ) $ such that $\X $ fits into an adjunction $L \dashv R \from \X \to \A $ satisfying the adhesivity condition (Def.~\ref{df:conditions}) and has a monic counit. We also require the grammar to induce the same language as the \emph{discrete grammar} obtained by replacing each rule $\ell \gets k \to r$ in $P$ with $\ell \gets LRk \to r$. It is already known that the category of directed graphs satisfy these conditions \cite[Prop.~3.3]{Ehrig_GraphGram}. We show this is also true of typed graphs, different flavors of hypergraphs, and Petri Nets.

Before proceeding, we reflect the above emphasis that this section is not about rewriting structured cospans. Instead, we will use structured cospans as a tool to inductively construct the language for a certain class of grammars.

\subsection{Discrete Grammars}
\label{sec:discrete-grammars}

In general, for any rewrite rule $ \ell \gets k \to r $, there is only one constraint on the value of $ k $: it must be a subobject of $ \ell $ and $ r $.  But requiring that $ k $ also be discrete can simplify any analysis involving that rule. This leads us question whether we can learn about a grammar $ ( \X,P ) $ by instead studying the grammar $ ( \X, P_\flat ) $, where $ P_\flat $ is obtained by \emph{discretizing} the apexes of every rule in $ P $. In this section we explain the $ \flat $ notation and make precise the concept of discreteness before giving the main result of this section that characterizes when $ ( \X,P ) $ and $ ( \X, P_\flat ) $ give the same language.  This result generalizes the characterization of discrete graph grammars given by Ehrig, et.~al. \cite[Prop.~3.3]{Ehrig_GraphGram}.

As a brief aside, experts in topos theory will know that discreteness comes from the flat modality on a local topos \cite[Ch.~3.6]{Johnstone_Sketches}).  While this is not unrelated to our definition of discreteness, we are really borrowing this term to use in this paper because it evokes the primary example of discretizing open graphs (Example \ref{ex:open-graph}).

\begin{definition}[Discrete comonad] \label{def:discrete-comonad}
  A comonad of an adjunction is called \df{discrete} if its counit is monic. 
\end{definition}

We can interpret a discrete comonad as returning largest interface $ \flat x $ supported by a network $ x $.  Here is an example illustrating how the discrete graph adjunction gives rise to $ \flat $.

\begin{example} \label{ex:open-graph}
  Consider the geometric morphism $$\adjunction{\Graph}{L}{R}{\Set}$$ defined by setting $ La $ to be the discrete graph on $ a $ and $ Rx $ to be the set of nodes in $ x $.  This adjunction induces the comonad $ \flat \bydef LR $ on $ \Graph $.  Applying $ \flat $ to a graph $ x $ returns the discrete graph underlying $ x $, for instance
  
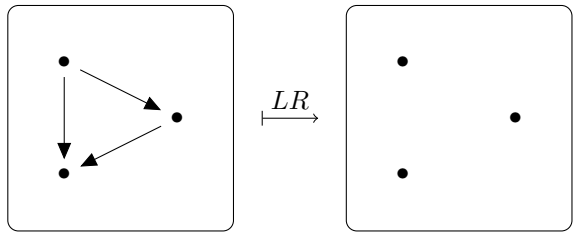
\[
    \begin{tikzpicture}[scale=0.5]
      \begin{scope}
        \node (a) at (0,2) {$ \bullet $};
        \node (b) at (2,1) {$ \bullet $};
        \node (c) at (0,0) {$ \bullet $};
        \draw [graph] 
        (a) edge[] (b)
        (a) edge[] (c)
        (b) edge[] (c);
        \draw [rounded corners] (-1,-1) rectangle (3,3);
      \end{scope}
      \begin{scope}[shift={(6,0)}]
        \node (a) at (0,2) {$ \bullet $};
        \node (b) at (2,1) {$ \bullet $};
        \node (c) at (0,0) {$ \bullet $};
        \draw [rounded corners] (-1,-1) rectangle (3,3);
      \end{scope}
      \draw [|->] (3.5,1) to node[above]{$ \flat $} (4.5,1);
    \end{tikzpicture}
  \]
  The counit $ \epsilon_x \from \flat x \to x $ is certainly monic as it includes the discrete graph $ \flat x $ into the graph $ x $, as in
  
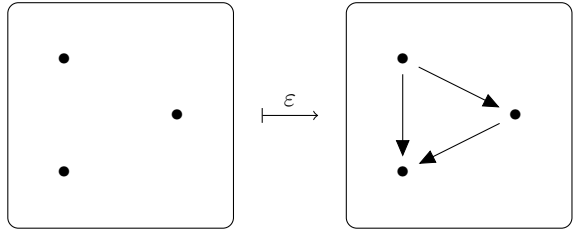
\[
    \begin{tikzpicture}[scale=0.5]
      \begin{scope}
        \node (a) at (0,2) {$ \bullet $};
        \node (b) at (2,1) {$ \bullet $};
        \node (c) at (0,0) {$ \bullet $};
        \draw [rounded corners] (-1,-1) rectangle (3,3);
      \end{scope}
      \begin{scope}[shift={(6,0)}]
        \node (a) at (0,2) {$ \bullet $};
        \node (b) at (2,1) {$ \bullet $};
        \node (c) at (0,0) {$ \bullet $};
        \draw [graph] 
        (a) edge[] (b)
        (a) edge[] (c)
        (b) edge[] (c);
        \draw [rounded corners] (-1,-1) rectangle (3,3);
      \end{scope}
      \draw [|->]
      (3.5,1) to node[above]{$ \epsilon $} (4.5,1);
    \end{tikzpicture}
  \]
\end{example}

Discrete comonads provide a tool to control the form of a grammar by replacing every rule $ \ell \gets k \to r $ with $ \ell \gets \flat k \to r $.

\begin{definition}[Discrete grammar] \label{def:DiscreteGrammar}
  Let $ \flat \from \X \to \X $ be a discrete comonad with counit  $ \epsilon $.  Given a grammar $ ( \X , P ) $, define $ P_\flat $ as the set containing \[ \ell \gets k \xgets{\epsilon} \flat k \xto{\epsilon} k \to r \] for each rule $ \ell \gets k \to r $ in $ P $. We call $ ( \X , P_\flat ) $ the \df{discrete grammar} underlying $ ( \X, P ) $.
\end{definition}

\subsection{Expressiveness of Grammars}
\label{sec:grammar-and-discrete}

The following lemma uses pushout complements to characterize when $(\X, P)$ and $(\X,P_\flat)$ generate the same language. We then use this characterization to show that this holds for both directed and undirected hypergraphs, Petri nets, and marked Petri nets. It also holds for certain slice categories, which allows us to extend this result to typed graphs and hypergraphs.

\begin{lemma} \label{thm:po-comp-same-rewrite-rel}
  Fix an adjunction \[\adjunction{\X}{L}{R}{\A}\] between adhesive categories with a monic counit and $L$ preserving pullbacks.  Take a grammar $(\X, P)$ and its underlying discrete grammar $(\X, P_{\flat})$. If for every derived rewrite rule
  \[
    \begin{tikzpicture}[
      baseline=(current bounding box.center)]
      scale=0.75]
      \node (l) at (0,1) {$ \ell $};
      \node (k) at (1,1) {$ k $};
      \node (r) at (2,1) {$ r $};
      \node (lprime) at (0,0) {$ \ell' $};
      \node (kprime) at (1,0) {$ k' $};
      \node (rprime) at (2,0) {$ r' $};
      \draw [cd]
      (k) edge[] (l)
      (k) edge[] (r)
      (kprime) edge[] (lprime)
      (kprime) edge[] (rprime)
      (l) edge[] (lprime)
      (k) edge[] (kprime)
      (r) edge[] (rprime);
      \draw (0.3,0.4) -- (0.4,0.4) -- (0.4,0.3);
      \draw (1.7,0.4) -- (1.6,0.4) -- (1.6,0.3);
    \end{tikzpicture}
  \]
  there is a pushout complement $d$ to the composite $\flat k \to k \to k'$, then $(\X,P)$ and $(\X, P_{\flat})$ generate the same language.
\end{lemma}

\begin{proof}
  Suppose that $\dderiv{\ell'}{r'}$ in the rewrite relation for $(\X,P)$. That means there is a double pushout diagram
  \[
    \begin{tikzpicture}[
      baseline=(current bounding box.center)]
      scale=0.75]
      \node (l) at (0,1) {$ \ell $};
      \node (k) at (1,1) {$ k $};
      \node (r) at (2,1) {$ r $};
      \node (lprime) at (0,0) {$ \ell' $};
      \node (kprime) at (1,0) {$ k' $};
      \node (rprime) at (2,0) {$ r' $};
      \draw [cd]
      (k) edge[] (l)
      (k) edge[] (r)
      (kprime) edge[] (lprime)
      (kprime) edge[] (rprime)
      (l) edge[] (lprime)
      (k) edge[] (kprime)
      (r) edge[] (rprime);
      \draw (0.3,0.4) -- (0.4,0.4) -- (0.4,0.3);
      \draw (1.7,0.4) -- (1.6,0.4) -- (1.6,0.3);
    \end{tikzpicture}
  \]
  in $\X$ so that the top span is in $P$. Then $\dderiv{\ell'}{r'}$ in $(\X, P_{\flat})$ as witnessed by the double pushout diagram
  \[
    \begin{tikzpicture}[
      baseline=(current bounding box.center)]
      scale=0.75]
      \node (l) at (0,1) {$ \ell $};
      \node (kleft) at (1,1) {$ k $};
      \node (kflat) at (2,1) {$ \flat k $};
      \node (kright) at (3,1) {$ k $};
      \node (r) at (4,1) {$ r $};
      \node (lprime) at (0,0) {$ \ell' $};
      \node (kprimeleft) at (1,0) {$ k' $};
      \node (d) at (2,0) {$ d $};
      \node (kprimeright) at (3,0) {$ k' $};
      \node (rprime) at (4,0) {$ r' $};
      \draw [cd]
      (kflat) edge[] (kleft)
      (kleft) edge[] (l)
      (kflat) edge[] (kright)
      (kright) edge[] (r)
      (d) edge[] (kprimeleft)
      (kprimeleft) edge[] (lprime)
      (d) edge[] (kprimeright)
      (kprimeright) edge[] (rprime)
      (l) edge[] (lprime)
      (kleft) edge[] (kprimeleft)
      (kflat) edge[] (d)
      (kright) edge[] (kprimeright)
      (r) edge[] (rprime);
      \draw (0.3,0.4) -- (0.4,0.4) -- (0.4,0.3);
      \draw (1.3,0.4) -- (1.4,0.4) -- (1.4,0.3);
      \draw (2.7,0.4) -- (2.6,0.4) -- (2.6,0.3);
      \draw (3.7,0.4) -- (3.6,0.4) -- (3.6,0.3);
    \end{tikzpicture}
  \] 

  Conversely, suppose that $\dderiv{\ell'}{r'}$ is generated by a rule $\spn{\ell}{\flat k}{r}$ in $P_\flat$.  That means there is a diagram
  \[
    \begin{tikzpicture}
      \node (l') at (0,0) {$ \ell' $};
      \node (l) at (0,1) {$ \ell $};
      \node (d) at (1,0) {$ d $};
      \node (flat k) at (1,1) {$ \flat k $};
      \node (r') at (2,0) {$ r' $};
      \node (r) at (2,1) {$ r $};
      \draw [->]
      (d) edge[] (l')
      (d) edge[] (r')
      (flat k) edge[] (l)
      (flat k) edge[] (r)
      (l) edge[] (l')
      (flat k) edge[] (d)
      (r) edge[] (r');
      \draw (0.3,0.4) -- (0.4,0.4) -- (0.4,0.3);
      \draw (1.7,0.4) -- (1.6,0.4) -- (1.6,0.3);
    \end{tikzpicture}
  \]
  By the construction of $P_\flat$, the maps $\flat k \to \ell$ and $\flat k \to r$ both factor through the counit $\flat k \to k$ and, consequently we get
  \[
    \begin{tikzpicture}
      \node (l') at (0,0) {$ \ell' $};
      \node (l) at (0,1) {$ \ell $};
      \node (k1) at (1,1) {$ k $};
      \node (d) at (2,0) {$ d $};
      \node (flat k) at (2,1) {$ \flat k $};
      \node (k2) at (3,1) {$ k $};
      \node (r') at (4,0) {$ r' $};
      \node (r) at (4,1) {$ r $};
      \draw [->]
      (d) edge[] (l')
      (d) edge[] (r')
      (flat k) edge[] (k1)
      (flat k) edge[] (k2)
      (k1) edge[] (l)
      (k2) edge[] (r) 
      (l) edge[] (l')
      (flat k) edge[] (d)
      (r) edge[] (r');
      \draw (0.3,0.4) -- (0.4,0.4) -- (0.4,0.3);
      \draw (3.7,0.4) -- (3.6,0.4) -- (3.6,0.3);
    \end{tikzpicture}
  \]
  Both maps $d \to \ell'$ and $d \to r'$ factor through $k +_{\flat k} d$ which results in the diagram
  \[
    \begin{tikzpicture}
      \node (l') at (0,0) {$ \ell' $};
      \node (l) at (0,1) {$ \ell $};
      \node (po) at (1.5,0) {$ k +_{\flat k} d$};
      \node (k) at (1.5,1) {$ k $};
      \node (r') at (3,0) {$ r' $};
      \node (r) at (3,1) {$ r $};
      \draw [->]
      (po) edge[] (l')
      (po) edge[] (r')
      (k) edge[] (l)
      (k) edge[] (r)
      (l) edge[] (l')
      (k) edge[] (po)
      (r) edge[] (r');
      \draw (0.3,0.4) -- (0.4,0.4) -- (0.4,0.3);
      \draw (2.7,0.4) -- (2.6,0.4) -- (2.6,0.3);
    \end{tikzpicture}
  \]
  whose pushouts result from the pushout pasting law. This witnesses $\dderiv{\ell'}{r'}$ derived from $P$.
\end{proof}

Showing that a grammar and its underlying discrete grammar generate the same language was initially shown for directed graphs.

\begin{example}[{\cite[Prop.~3.3]{Ehrig_GraphGram}}]
  Returning to the discrete graph geometric morphism $$\adjunction{\Graph}{L}{R}{\Set},$$ any grammar $(\Graph, P)$ generates the same language as the discrete grammar $(\Graph, P_{\flat})$ .
\end{example}

\subsubsection{Hypergraphs}

We can now show for the first time that this result also holds for hypergraphs, both directed and undirected. Following Bonchi, et.~al. \cite{bonchi-etal_graph-rewriting-interfaces}, we realize hypergraphs as a presheaf topos in the following way.

\begin{definition}
  Define the category of directed hypergraphs $\cat{DHGraph}$ as the presheaf category $[\cat{DH}\op, \Set]$ where $\cat{DH}$ has objects $\NN \times \NN + \{N\}$ and arrows \[\cat{DH}(N,(n,m)) \bydef \{s_1, \dotsc, s_n\} + \{t_1, \dotsc, t_m\}.\] All other homsets are trivial. For such a presheaf, the image of $(n,m) \in \NN \times \NN$ gives the set of edges with in-degree $n$ and  out-degree $m$, the image of $N$ gives the set of nodes, the image of each \emph{source} map $s_j^{\op} \from (n,m) \to N$  records the $j$-th incidence for the source of each edge, and the  image of each \emph{target} map $t_j^{\op} \from (n,m) \to N$  records the $j$-th incidence for the target of each edge.
  
  Define the category of undirected hypergraphs $\cat{HGraph}$ as the category of presheaves $[\cat{H}\op,\Set]$ on the category $\cat{H}$ with objects $\NN + \{N\}$ and arrows \[\cat{H}(N,n) \bydef \{x_1, \dotsc, x_n\}.\] All other homsets are trivial. For such a presheaf, the image of $n \in \NN$ gives the set of edges with degree $n$, the image of $N$ gives the set of nodes, and the image of each map $x_j^{\op} \from n \to N$ records the $j$-th incidence for each edge.
\end{definition}

Both undirected hypergraphs can be \emph{discretized} via geometric morphisms $$\adjunction{\DHGraph}{L}{R}{\Set} \quad\quad \adjunction{\HGraph}{L}{R}{\Set}$$ where $R$ returns the node set and $L$ returns a discrete hypergraph on a given set. It is straightforward to check that, in each case, the counit of the induced comonad is indeed monic. Therefore, for a grammar of directed hypergraphs $(\DHGraph, P)$, we have the underlying discrete grammar $(\DHGraph, P_{\flat})$. This holds too for undirected hypergraphs.

\begin{proposition} \label{prop:dir-hyps-discrete-gram}
  A grammar of directed hypergraphs $(\DHGraph, P)$ and its underlying discrete grammar $(\DHGraph, P_{\flat})$ have the same language.
\end{proposition}

\begin{proof}
  Consider a derived rewrite rule 
  \[
    \begin{tikzpicture}[
      baseline=(current bounding box.center)]
      scale=0.75]
      \node (1) at (0,1) {$ \ell $};
      \node (2) at (1,1) {$ k $};
      \node (3) at (2,1) {$ r $};
      \node (4) at (0,0) {$ \ell' $};
      \node (5) at (1,0) {$ k' $};
      \node (6) at (2,0) {$ r' $};
      \draw [cd]
      (2) edge (1)
      (2) edge (3)
      (5) edge (4)
      (5) edge (6)
      (1) edge node[left]{$ m $} (4)
      (2) edge (5)
      (3) edge (6); 
      \draw (0.3,0.4) -- (0.4,0.4) -- (0.4,0.3);
      \draw (1.7,0.4) -- (1.6,0.4) -- (1.6,0.3);
    \end{tikzpicture}
  \]
By Lemma \ref{thm:po-comp-same-rewrite-rel}, it suffices to find a hypergraph $d$ that is a pushout complement to $\flat k \to k \to k'$. Define $d$ by $d(N) \bydef k'(N)$ and $d(n,m) \bydef k'(n,m) - k(n,m)$ for each in-degree, out-degree pair $(n,m)$. Recalling that the directed hypergraph $\flat k$ has no edges and the same nodes as $k$, we have pushouts
  \[
    \begin{tikzpicture}
      \node (k) at (0,1) {$ k(n,m) $};
      \node (k') at (0,0) {$ k'(n,m) $};
      \node (LRk) at (2,1) {$ \flat k(n,m) $};
      \node (d) at (2,0) {$ d(n,m) $};
      \draw [->]
      (LRk) edge[] (k)
      (LRk) edge[] (d)
      (k) edge[] (k')
      (d) edge[] (k');
      \draw (0.3,0.4) -- (0.4,0.4) -- (0.4,0.3);
    \end{tikzpicture}
    \quad\quad\quad
    \begin{tikzpicture}
      \node (k) at (0,1) {$ k(N) $};
      \node (k') at (0,0) {$ k'(N) $};
      \node (LRk) at (2,1) {$ \flat k(N) $};
      \node (d) at (2,0) {$ d(N) $};
      \draw [->]
      (LRk) edge[] (k)
      (LRk) edge[] (d)
      (k) edge[] (k')
      (d) edge[] (k');
      \draw (0.3,0.4) -- (0.4,0.4) -- (0.4,0.3);
    \end{tikzpicture}
  \]
  ensuring that $d$ is the pushout complement we sought.
\end{proof}

\begin{proposition} \label{prop:undir-hyps-discrete-gram}
  A grammar of undirected hypergraphs $(\HGraph, P)$ and its underlying discrete grammar $(\HGraph, P_{\flat})$ have the same language.
\end{proposition}

\begin{proof}
  The proof follows the same reasoning as the directed hypergraph case.
\end{proof}

\subsubsection{Petri Nets, Unmarked and Marked}

Petri nets have various definitions in the literature, particularly their morphisms. Some authors define morphisms to capture the \emph{behavior} of Petri Nets \cite{baez-master_open-pnets}.  Our interests lie in the structure of Petri Nets as opposed to the behavior, so we take the following definition as used by Johnstone, et.~al. \cite{johnstone_quasi}.

\begin{definition}[Petri Nets]
  A \df{Petri Net} is a tuple $N \bydef (P,T,s,t)$ where $P$ is a set of places, $T$ is a set of \emph{transitions}, and $s,t \from T \to P^\oplus$ are maps that encode the \emph{source} and \emph{target} of the transitions by assigning each transition to a collection places. A \df{morphism of Petri Nets} is a pair $(f_P,f_T) \from N \to N'$ where $f_P \from P \to P'$ and $f_T \from T \to T'$ preserve sources and targets: $s'f_T=f_{P}^{\oplus}s$ and $t'f_T=f_{P}^{\oplus}t$.  The category of Petri Nets and their morphisms is denoted by $\cat{PN}$. This category is equivalent to the category of directed hypergraphs.
\end{definition}

\begin{proposition} 
  A grammar of Petri Nets $(\cat{PN}, P)$ and its underlying discrete grammar $(\cat{PN}, P_{\flat})$ have the same language.
\end{proposition}

\begin{proof}
   The result follows from a similar argument as in Prop.~\ref{prop:dir-hyps-discrete-gram}.
\end{proof}

\begin{definition}[Marked Petri Nets]
  A \df{marked Petri net} is a tuple $M \bydef (N,K,k)$ where $N \bydef (P,T,s,t)$ is a Petri net, $K$ is a set of \emph{tokens}, and $k \from K \to P$ is a \emph{marking} of the places by the tokens. A \df{morphism of marked Petri nets} $(f_N,f_K) \from M \to M'$ is a Petri net morphism $f_N \from N \to N'$ together with a map between token sets $f_K \from K \to K'$ such that $k'f_K = f_Pk$. The category of marked Petri nets and their morphisms is denoted by $\cat{MPN}$. This is equivalent to a presheaf category $[\cat{M}\op,\Set]$ for $\cat{M}$ the category with objects $\NN \times \NN + \{P,K\}$ and arrows $\cat{M}(P,(n,m)) = \{s_1, \dotsc, s_n\} + \{t_1, \dotsc, t_m\}$ plus $\cat{M}(P,K) = \{k\}$.
\end{definition}

A \emph{discrete marked Petri net} is one without transitions or tokens. Similar to graphs and hypergraphs, there is a discrete marked Petri net adjunction. 

\begin{proposition} 
  A grammar of marked Petri Nets $(\cat{MPN}, P)$ and its underlying discrete grammar $(\cat{MPN}, P_{\flat})$ have the same language.
\end{proposition}

\begin{proof}
  The result follows from a similar argument as in Prop.~\ref{prop:dir-hyps-discrete-gram}.
\end{proof}

\subsubsection{Typed Objects}

Occasionally, those working with rewriting \emph{typify} their objects. For instance, instead of working with graphs, one might work with graphs whose edges and notes each come in multiple colors. This can be accomplished by working in an appropriate slice category.  For example, suppose we wanted to work with graphs whose edges come in red $r$ and green $g$ colors. Add this information to a graph using a map to
\[
  \begin{tikzpicture}
    \node (b) at (0,0) {$ \bullet $};  
    \draw [->,loop above,looseness=15] (b) edge node[left]{$ r $} (b);
    \draw [->,loop below,looseness=15] (b) edge node[left]{$ g $} (b);
  \end{tikzpicture}
\]
so that the images of the edges and nodes of $g$ encode the typing information.  In short, we may want to be rewriting in a slice category $\X / x$ using an adhesive category $\X$ and object $x$.

This next result provides conditions for a grammar and its underlying discrete grammar in $\X/x$ to generate the same language.  The first condition is that the ``untyped grammars'' $(\X, P)$ and $(\X,P_\flat)$ generate the same rewrite relation.  To move to a grammar in $\X/x$, there would be user-driven decisions to make, namely the how the rewrite rules in $P$ will be mapped to $x$.

\begin{proposition}
  Fix an adjunction $L \dashv R \from \X \to \A$ between adhesive categories with a monic counit and $L$ preserving pullbacks. Fix an object $x$ in $\X$. Consider both the comonad $\flat \coloneqq LR$ and the induced adjunction between slice categories $L_x \dashv R_x \from \X/x \to \A/Rx$ that also has a monic counit and $L_x$ preserving pullbacks. Let $(\X/x, P_x)$ be any grammar, $(\X/x, P_{x,\flat})$ be the underlying discrete grammar, and  $(\X, P)$ and $(\X, P_{\flat})$ be their projections into $\X$.

  If for every derived rewrite rule
  \[
    \begin{tikzpicture}[
      baseline=(current bounding box.center)]
      scale=0.75]
      \node (l) at (0,1) {$ \ell $};
      \node (k) at (1,1) {$ k $};
      \node (r) at (2,1) {$ r $};
      \node (lprime) at (0,0) {$ \ell' $};
      \node (kprime) at (1,0) {$ k' $};
      \node (rprime) at (2,0) {$ r' $};
      \draw [cd]
      (k) edge[] (l)
      (k) edge[] (r)
      (kprime) edge[] (lprime)
      (kprime) edge[] (rprime)
      (l) edge[] (lprime)
      (k) edge[] (kprime)
      (r) edge[] (rprime);
      \draw (0.3,0.4) -- (0.4,0.4) -- (0.4,0.3);
      \draw (1.7,0.4) -- (1.6,0.4) -- (1.6,0.3);
    \end{tikzpicture}
  \]
  in $(\X,P)$ there is a pushout complement $d$ to the composite $\flat k \to k \to k'$, then $(\X/x,P_x)$ and $(\X/x, P_{x,\flat})$ generate the same rewrite relation.
\end{proposition}

\begin{proof}
  Recalling that pushouts in $\X/x$ are computing by projecting to $\X$, this result follows from a similar argument as in Lemma \ref{thm:po-comp-same-rewrite-rel}.
\end{proof}

\begin{corollary}
  Grammars and their underlying discrete grammars generate the same rewrite relation for typed directed graphs, typed directed hypergraphs, and typed undirected hypergraphs.  
\end{corollary}

\subsection{Inductive Rewriting Systems}
\label{sec:inductive-rewriting-hypergraphs}

In the networks perspective, a grammar $( \X,P )$ can be thought of as comprising an adhesive category $\X$ of closed networks and a set of rules $P$ stating how to revise the closed networks. For example, one might revise an network of resistors by replacing a series of resistors with a single resistor. An important part of our construction involves decomposing a closed network $x$ into other networks $x_1, \dotsc, x_n$ that are somehow connected. Structured cospans provide a way to form these connections, hence we want $\X$ to fit into an appropriate adjunction $L \dashv R$. Because our construction uses structured cospans, we need a way to represent a closed network with structured cospans. We turn a closed network $x$ into the structured cospan $\csp{L0}{x}{L0}$ with an empty interface. Note that here $0$ denotes an initial object in $\A$, and $L0$ is an initial object in $\X$ because $L$ is a left adjoint.

The particular decomposition we use is determined by the grammar. That is, we start with a grammar $( \X,P )$ where $\X$ can fit into an adjunction $L \dashv R \from \X \to \A$ with a monic counit satisfying the adhesivity condition. This gives a discrete comonad $\flat \bydef LR$ which allows us to form the discrete grammar $(\X,P_\flat)$ as in Definition \ref{def:DiscreteGrammar}.

Finally, we need to define a particular grammar used to generate our inductive viewpoint of rewriting.

\begin{definition} \label{df:decompGrammar}
  Fix an adjunction $L \dashv R \from \X \to \A$ and a grammar $G \bydef (\X, P)$. Define $ G' = ( _L \StrCsp , P'_\flat ) $ to be a structured cospan grammar where $ P'_\flat $ contains the rule
  \[
     \begin{tikzpicture}[baseline=(current  bounding  box.center)]
       \node (2) at (1,2) {$ \ell $};
       \node (4) at (0,1) {$ L 0 $};
       \node (5) at (1,1) {$ \flat k $};
       \node (6) at (2,1) {$ \flat k $};
       \node (8) at (1,0) {$ r $};
       \draw [cd] (4) to (5);
       \draw [cd] (6) to (5);
       \draw [cd] (4) to (2);
       \draw [cd] (4) to (8);
       \draw [cd] (5) to (2);
       \draw [cd] (5) to (8);
       \draw [cd] (6) to (2);
       \draw [cd] (6) to (8);
     \end{tikzpicture}
   \]
  for each rule $ \spn{\ell}{\flat k}{r} $ of $ P_{\flat} $. 
\end{definition}

Applying the language construction $\LLang (G')$, then, encodes $\Lang (G)$ in its 2-cells which are generated by the rules in $P_\flat$. Since $\Lang (G)$ is generated by the rewrite relation for $G$, this forms the inductive viewpoint of the rewrite relation.

\begin{therm} \label{thm:inductive-rewriting}
  Let $ L \dashv R \from \X \to \A $ be an adjunction that satisfies the adhesivity condition and induces a discrete comonad $\flat$. Given any grammar $G \bydef (\X,P)$ whose language is the same as its underlying discrete grammar $G_\flat \bydef (\X,P_\flat)$. There is an equivalence between $\Lang (G)$ and $\LLang (G')(L0,L0)$
\end{therm}

\begin{proof}
  We show sufficiency by inducting on the length of the generating arrows of $\Lang (G)$, which come from the rewrite relation $\deriv{}{}$. If $ \deriv{g}{h} $ in a single step, meaning that there is a diagram
  \[
    \begin{tikzpicture}
      \node (1t) at (0,1) {$ \ell $};
      \node (2t) at (1,1) {$ \flat k $};
      \node (3t) at (2,1) {$ r $};
      \node (1b) at (0,0) {$ g $};
      \node (2b) at (1,0) {$ d $};
      \node (3b) at (2,0) {$ h $};
      \draw [cd] (2t) to (1t);
      \draw [cd] (2t) to (3t);
      \draw [cd] (2b) to (1b);
      \draw [cd] (2b) to (3b);
      \draw [cd] (1t) to (1b);
      \draw [cd] (2t) to (2b);
      \draw [cd] (3t) to (3b);
      \draw (0.3,0.4) -- (0.4,0.4) -- (0.4,0.3);
      \draw (1.7,0.4) -- (1.6,0.4) -- (1.6,0.3);
    \end{tikzpicture}
  \]
  in $\X$ whose top row is a rewrite rule in $P$, then the desired 2-cell is the horizontal composition of
  \[
    \begin{tikzpicture}[baseline =(current bounding box.center)]
      \node (2t) at (1,2) {$ \ell $};
      \node (4t) at (3,2) {$ d $};
      \node (1m) at (0,1) {$ L0 $};
      \node (3m) at (2,1) {$ LRk $};
      \node (5m) at (4,1) {$ L0 $};
      \node (2b) at (1,0) {$ r $};
      \node (4b) at (3,0) {$ d $};
      \draw [cd] (1m) to (2t);
      \draw [cd] (1m) to (2b);
      \draw [cd] (3m) to (2t);
      \draw [cd] (3m) to (2b);
      \draw [cd] (3m) to (4t);
      \draw [cd] (3m) to node [] {\scriptsize{$  $}} (4b);
      \draw [cd] (5m) to node [] {\scriptsize{$  $}} (4t);
      \draw [cd] (5m) to node [] {\scriptsize{$  $}} (4b);
      \draw[double,->] (1,1.25) -- node [right] {\scriptsize{$ \gamma_p $}}  (1,0.75);
      \draw[double,->] (3,1.25) -- node [right] {\scriptsize{$ \id $}}(3,0.75);
    \end{tikzpicture}
    \mapsto
    \begin{tikzpicture}[baseline =(current bounding box.center)]    
      \node (l) at (0,1) {$ L 0 $};
      \node (r) at (2,1) {$ L 0 $};
      \node (b) at (1,0) {$ r +_{LRk} d $};
      \node (t) at (1,2) {$ \ell +_{LRk} d $};
      \draw [cd]
      (l) to (t)
      (l) to (b)
      (r) to (t)
      (r) to (b);    
      \draw[double,->] (1,1.25) -- node [right] {\scriptsize{$ \gamma_p \id $}}  (1,0.75);   
    \end{tikzpicture}
  \]
  where $g \cong l +_{\flat k} d$ and $h \cong r +_{\flat k} d$.  Then the 2-cell for a derivation $ \dderiv{\deriv{g}{h}}{j} $ is the vertical composition of
  
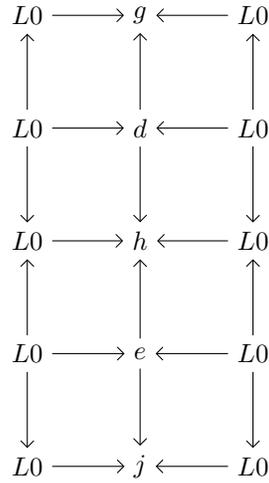
\[
    \begin{tikzpicture}
      \node (L0l) at (0,1) {$ L 0 $};
      \node (L0r) at (2,1) {$ L 0 $};
      \node (g) at (1,2) {$ g $};
      \node (h) at (1,1) {$ h $};
      \node (j) at (1,0) {$ j $};
      \draw [cd] (L0l) to (g);
      \draw [cd] (L0l) to (h);
      \draw [cd] (L0l) to (j);
      \draw [cd] (L0r) to (g);
      \draw [cd] (L0r) to (h);
      \draw [cd] (L0r) to (j);
      \draw [double,->] (1,1.65) -- (1,1.35);
      \draw [double,->] (1,0.65) -- (1,0.35);
    \end{tikzpicture}
  \]
  The top square is from $ \deriv{g}{h} $ and the second from $ \dderiv{h}{j} $.

  Conversely, we proceed by structural induction to show that given an arrow in the hom-category $ \LLang (G')(L0,L0) $, there is a corresponding derivation per the rewrite relation. The base case is the generating 2-cells which holds by construction. Next, we show that both horizontal and vertical composition preserve $\deriv{}{}$. Vertical composition preserves $\deriv{}{}$ due to the transitivity of the rewrite relation. Horizontal composition must have the form
  \[
    \LLang (G')(L0,La) \times \LLang (G') (La,L0) \to \LLang (G')(L0,L0).
  \]
  So given
  \[
    \begin{tikzpicture}[baseline =(current bounding box.center)]
      \node (2t) at (1,2) {$ g $};
      \node (4t) at (3,2) {$ g' $};
      \node (1m) at (0,1) {$ L0 $};
      \node (3m) at (2,1) {$ La $};
      \node (5m) at (4,1) {$ L0 $};
      \node (2b) at (1,0) {$ h $};
      \node (4b) at (3,0) {$ h' $};
      \draw [cd] (1m) to (2t);
      \draw [cd] (1m) to (2b);
      \draw [cd] (3m) to (2t);
      \draw [cd] (3m) to (2b);
      \draw [cd] (3m) to (4t);
      \draw [cd] (3m) to node [] {\scriptsize{$  $}} (4b);
      \draw [cd] (5m) to node [] {\scriptsize{$  $}} (4t);
      \draw [cd] (5m) to node [] {\scriptsize{$  $}} (4b);
      \draw[double,->] (1,1.25) -- node [right] {\scriptsize{$ \gamma_0 $}}  (1,.75);
      \draw[double,->] (3,1.25) -- node [right] {\scriptsize{$ \gamma_1 $}}(3,.75);
    \end{tikzpicture}
  \]
  where $\deriv{g}{h}$ and $\deriv{g'}{h'}$ then we get  $\deriv{g+_{a}g'}{h+_{a} h'}$ by realizing the horizontal   composition as
  \[
    \begin{tikzpicture}[baseline =(current bounding box.center)]
      \node (2t) at (1,2) {$ g $};
      \node (4t) at (3,2) {$ g' $};
      \node (1m) at (0,1) {$ L0 $};
      \node (2m) at (1,1) {$ h $};
      \node (3m) at (2,1) {$ LRk $};
      \node (4m) at (3,1) {$ g' $};
      \node (5m) at (4,1) {$ L0 $};
      \node (2b) at (1,0) {$ h $};
      \node (4b) at (3,0) {$ h' $};
      \draw [cd] (1m) to (2t);
      \draw [cd] (1m) to (2m);
      \draw [cd] (1m) to (2b);
      \draw [cd] (3m) to (2t);
      \draw [cd] (3m) to (2m);
      \draw [cd] (3m) to (2b);
      \draw [cd] (3m) to (4t);
      \draw [cd] (3m) to (4m);
      \draw [cd] (3m) to (4b);
      \draw [cd] (5m) to (4t);
      \draw [cd] (5m) to (4m);
      \draw [cd] (5m) to (4b);
      \draw[double,->] (1,1.65) -- node [right] {\scriptsize{$ \gamma $}}  (1,1.35);
      \draw[double,->] (3,1.65) -- node [right] {\scriptsize{$ \id $}}  (3,1.35);
      \draw[double,->] (1,.65) -- node [right] {\scriptsize{$ \id $}} (1,0.35);
      \draw[double,->] (3,.65) -- node [right] {\scriptsize{$ \gamma' $}}(3,0.35);
    \end{tikzpicture}
  \]
  composing horizontally, then performing this rewrite by first   applying $\deriv{g}{h}$ to $g$ to get   $\deriv{g+_{La}g'}{h+_{La} g'}$ and then applying $\deriv{g'}{h'}$ to $g'$ to get $\deriv{g+_{La}g'}{h+_{La} h'}$.
\end{proof}

As a corollary, we have that Theorem \ref{thm:inductive-rewriting} applies to grammars in the categories of directed hypergraphs, undirected hypergraphs, unmarked Petri Nets, marked Petri nets, and their typed variants.  It also applies to typed graphs so applies to the ZX-calculus \cite{ZX}. 

\section{Conclusion}
\label{sec:conclusion}

We have introduced a new category of structured cospans $_L\StrCsp$ and characterized when it is adhesive and when it is a topos. Under these conditions, $_L\StrCsp$ admits a theory of rewriting. Next, we have provided a condition for a grammar and its underlying discrete grammar to induce the same language. Using that condition, we introduced an inductive viewpoint for rewriting, thus allowing for proof strategies involving structural induction. In particular, we have shown this holds for directed hypergraphs, undirected hypergraphs, their typed variants, for Petri nets, marked Petri nets, and the ZX-calculus.  

\bibliographystyle{plain}
\bibliography{assets/biblio}

\end{document}